\documentstyle[twoside]{article}
\title{A Poisson-Jacobi-type transformation for the sum $\sum_{n=1}^\infty n^{-2m}\exp (-an^2)$ for positive integer $m$}
\author{\sc R. B.\ Paris \\
{\em School of Engineering, Computing and Applied Mathematics}, \\
{\em University of Abertay Dundee, Dundee DD1 1HG, UK}
}
\begin{document}
\def\f#1#2{\mbox{${\textstyle \frac{#1}{#2}}$}}
\def\dfrac#1#2{\displaystyle{\frac{#1}{#2}}}
\def\boldal{\mbox{\boldmath $\alpha$}}
{\newcommand{\Sgoth}{S\;\!\!\!\!\!/}
\newcommand{\bee}{\begin{equation}}
\newcommand{\ee}{\end{equation}}
\newcommand{\lam}{\lambda}
\newcommand{\ka}{\kappa}
\newcommand{\al}{\alpha}
\newcommand{\th}{\theta}
\newcommand{\fr}{\frac{1}{2}}
\newcommand{\fs}{\f{1}{2}}
\newcommand{\g}{\Gamma}
\newcommand{\br}{\biggr}
\newcommand{\bl}{\biggl}
\newcommand{\ra}{\rightarrow}
\newcommand{\mbint}{\frac{1}{2\pi i}\int_{c-\infty i}^{c+\infty i}}
\newcommand{\mbcint}{\frac{1}{2\pi i}\int_C}
\newcommand{\mboint}{\frac{1}{2\pi i}\int_{-\infty i}^{\infty i}}
\newcommand{\gtwid}{\raisebox{-.8ex}{\mbox{$\stackrel{\textstyle >}{\sim}$}}}
\newcommand{\ltwid}{\raisebox{-.8ex}{\mbox{$\stackrel{\textstyle <}{\sim}$}}}
\renewcommand{\topfraction}{0.9}
\renewcommand{\bottomfraction}{0.9}
\renewcommand{\textfraction}{0.05}
\newcommand{\mcol}{\multicolumn}
\date{}
\maketitle
\pagestyle{myheadings}
\markboth{\hfill \sc R. B.\ Paris  \hfill}
{\hfill \sc  Poisson-Jacobi-type transformation\hfill}
\begin{abstract}
We obtain an asymptotic expansion for the sum
\[S(a;w)=\sum_{n=1}^\infty \frac{e^{-an^2}}{n^{w}}\]
as $a\rightarrow 0$ in $|\arg\,a|<\fs\pi$ for arbitrary finite $w>0$. The result when $w=2m$, where $m$ is a positive integer, is the analogue of the well-known Poisson-Jacobi transformation
for the sum with $m=0$.  Numerical results are given to illustrate the accuracy of the expansion.

\vspace{0.4cm}

\vspace{0.3cm}

\end{abstract}

\vspace{0.3cm}

\noindent $\,$\hrulefill $\,$

\vspace{0.2cm}

\begin{center}
{\bf 1. \  Introduction}
\end{center}
\setcounter{section}{1}
\setcounter{equation}{0}
\renewcommand{\theequation}{\arabic{section}.\arabic{equation}}
\vspace{0.6cm}
The classical Poisson-Jacobi transformation is given by
\bee\label{e11}
\sum_{n=1}^\infty e^{-an^2}=\frac{1}{2}\sqrt{\frac{\pi}{a}}-\frac{1}{2}+\sqrt{\frac{\pi}{a}} \sum_{n=1}^\infty e^{-\pi^2n^2/a},
\ee
where the parameter $a$ satisfies $\Re (a)>0$. This transformation relates a sum of Gaussian exponentials involving the parameter $a$ to a similar sum with parameter $\pi^2/a$. In the case $a\ra 0$ in $\Re (a)>0$, 
the convergence of the sum on the left-hand side becomes slow, whereas the sum on the right-hand side converges rapidly in this limit. Various proofs of the well-known result (\ref{e11}) exist in the literature; see, for example, 
\cite[p.~120]{PK}, \cite[p.~60]{T1} and  \cite[p.~124]{WW}.

In this note we consider the sum
\bee\label{e12}
S(a;w)=\sum_{n=1}^\infty \frac{e^{-an^2}}{n^w}\qquad (\Re (a)>0).
\ee
This sum converges for any finite value of the parameter $w$ provided $\Re (a)>0$; when $a=0$ then $S(0;w)$ reduces to the Riemann zeta function $\zeta(w)$ when $\Re (w)>1$. Consequently, the series in (\ref{e12}) can be viewed as a smoothed Dirichlet series for $\zeta(w)$. The asymptotic expansion of $S(a;w)$ as $a\ra 0$ in $\Re (a)>0$ is straightforward. The most interesting case arises when $w=2m$, where $m$ is a positive integer, for which we establish a transformation
for $S(a;2m)$ analogous to that in (\ref{e11}) valid as $a\ra 0$ in $\Re (a)>0$. This similarly involves the series in (\ref{e12}) with $a$ replaced by $\pi^2/a$, but with each term decorated by an asymptotic series in $a$. A recent application of the series with $w=2$ and $w=4$ has arisen in the geological problem of thermochronometry  in spherical geometry \cite{WFK}.
\vspace{0.6cm}

\begin{center}
{\bf 2. \ An expansion for $S(a;w)$ as $a\ra 0$ when $w\neq 2,4, \ldots$}
\end{center}
\setcounter{section}{2}
\setcounter{equation}{0}
\renewcommand{\theequation}{\arabic{section}.\arabic{equation}}
Our starting point is the well-known Cahen-Mellin integral (see, for example, \cite[\S 3.3.1]{PK})
\bee\label{e21}
z^{-\alpha} e^{-z}=\frac{1}{2\pi i}\int_{c-\infty i}^{c+\infty i} \g(s-\alpha) z^{-s}ds \qquad(z\neq 0,\ |\arg\,z|<\fs\pi),
\ee
where $c>\Re (\alpha)$ so that the integration path passes to the right of all the poles of $\g(s-\alpha)$ situated at $s+\alpha-k$ ($k=0, 1, 2, \ldots$).
For simplicity in presentation we shall assume throughout real values of $w>0$.
Then, it follows that
\begin{eqnarray*}
S(a;w)&=&\sum_{n=1}^\infty \frac{e^{-an^2}}{n^w}
=\sum_{n=1}^\infty\frac{n^{-w}}{2\pi i}\int_{c-\infty i}^{c+\infty i} \g(s) (an^2)^{-s}ds\\
&=&\frac{1}{2\pi i}\int_{c-\infty i}^{c+\infty i} \g(s) \zeta(2s+w)a^{-s}ds,
\end{eqnarray*}
upon reversal of the order of summation and integration, which is justified when $c>\max\{0,\fs-\fs w\}$, and evaluation of the inner sum in terms of the Riemann zeta function.
The integrand possesses simple poles at $s=\fs-\fs w$ and $s=-k$ ($k=0, 1, 2, \ldots$), except if $w=2m+1$ is an odd positive integer when the pole at $s=\fs-\fs w$ is double. The case when $w=2m$ is an even positive integer requires
a separate investigation which is discussed in Section 3.

Consider the integral taken round the rectangular contour with vertices at $c\pm iT$, $-c'\pm iT$, where $c'>0$. The contribution from the upper and lower sides $s=\sigma\pm iT$, $-c'\leq\sigma\leq c$, vanishes as $T\ra\infty$ provided $|\arg\,a|<\fs\pi$, since from the behaviour
\[\g(\sigma\pm it)=O(t^{\sigma-\fr}e^{-\fr\pi t}),\qquad \zeta(\sigma\pm it)=O(t^{\mu(\sigma)}\log^A t),\qquad(t\ra\infty),\]
where for $\sigma$ and $t$ real
\[\mu(\sigma)=0\ (\sigma>1),\quad \fs-\fs\sigma\ (0\leq\sigma\leq 1),\quad \fs-\sigma\ (\sigma<0),\]
\[A=1\ (0\leq\sigma\leq 1),\quad A=0\ \mbox{otherwise},\]
the modulus of the integrand is controlled by $O(T^{\sigma+\mu(\sigma)-\fr}\log\,T e^{-\Delta T})$, with $\Delta=\fs\pi-|\arg\,a|$. The residue at the double pole $s=-m$ when $w=2m+1$ ($m=0, 1, 2,\ldots$) is given by
\[\frac{(-a)^m}{m!} \{\gamma-\fs\log\,a+\fs\psi(m+1)\},\]
where $\gamma$ is Euler's constant and $\psi(x)$ is the logarithmic derivative of the gamma function. Displacement of the integration path to the left over the poles then yields (provided $w\neq 2m$)
\bee\label{e22}
S(a;w)=J(a;w)+\sum_{k=0}^{N-1}{}^{\!'}\frac{(-)^k}{k!}\,\zeta(w-2k)a^k+R_N,
\ee
where 
\[J(a;w)=\left\{\begin{array}{ll} \fs\g(\fs-\fs w) a^{(w-1)/2} & (w\neq 2m+1) \\
\\
\dfrac{(-a)^m}{m!} \{\gamma-\fs\log\,a+\fs\psi(m+1)\} & (w=2m+1), \end{array}\right.,\]
$N$ is a positive integer such that $N>\fs w+\fs$ and the prime on the sum over $k$ denotes the omission of the term corresponding to $k=m$ when $w=2m+1$.

The remainder $R_N$ is
\bee\label{e22a}
R_N=\frac{1}{2\pi i}\int_{-c-\infty i}^{-c+\infty i}\g(s) \zeta(w+2s) a^{-s}ds,\qquad c=N-\fs.
\ee
It is shown in the appendix, when $w\neq 2, 4, \ldots\,$, that $R_N=O(a^{N-\fr})$ as $a\ra 0$ in $|\arg\,a|<\fs\pi$, with the constant implied in the $O$-symbol growing at least like $\g(N\!+\!1\!-\!\fs w)$. This establishes that the above series over $k$ diverges as $N\ra\infty$ and that (\ref{e22}) is therefore an asymptotic expansion.

We remark that the algebraic expansion (\ref{e22}) also contains a subdominant exponentially small component as $a\ra 0$; compare \cite[\S 8.1.5]{PK} for the particular case $w=0$. We do not consider this further in the present paper.
\vspace{0.6cm}

\begin{center}
{\bf 3. \ An expansion for $S(a;2m)$ when $m=1, 2, \ldots$}
\end{center}
\setcounter{section}{3}
\setcounter{equation}{0}
\renewcommand{\theequation}{\arabic{section}.\arabic{equation}}
The case $w=2m$, where $m$ is a positive integer, is more interesting as this leads to the analogue of the Poisson-Jacobi transformation (\ref{e11}). There is now only a finite set of poles of the integrand in (\ref{e21})
at $s=\fs-\fs w$ and $s=0, -1, -2, \ldots , -m$, since the poles of $\g(s)$ at $s=-m-k$ ($k=1, 2, \ldots$) are cancelled by the trivial zeros of the zeta function $\zeta(2m+2s)$ at $s=-m-1, -m-2, \ldots\,$. This has the consequence that the integrand is holomorphic in $\Re (s)<-m$, so that further displacement of the contour can produce no additional algebraic terms in the expansion of $S(a;2m)$. Thus, we find when $w=2m$
\bee\label{e23}
S(a;2m)=\fs\g(\fs-m) a^{m-\fr}+\sum_{k=0}^m\frac{(-)^k}{k!}\,\zeta(2m-2k)\,a^k +I_L,
\ee
where, upon making the change of variable $s\ra -s$,
\bee\label{e24}
I_L=\frac{1}{2\pi i}\int_L \g(-s)\zeta(2m-2s)a^{s}ds
\ee
and $L$ denotes a path parallel to the imaginary axis with $\Re (s)>m$. 

We now employ the functional relation for $\zeta(s)$ given by \cite[p.~269]{WW}
\bee\label{e23a}
\zeta(s)=2^s\pi^{s-1} \zeta(1-s) \g(1-s) \sin \fs\pi s
\ee
to convert the argument of the zeta function in (\ref{e24}) into one with real part greater than unity. The integral in (\ref{e24}) can then be written in the form
\[\frac{(-)^m (2\pi)^{2m}}{2\pi i}\int_L \zeta(2s-2m+1)\,\frac{\g(2s-2m+1)}{\g(s+1)}\,\left(\frac{a}{4\pi^2}\right)^s ds.\]
Since on the integration path $\Re (2s-2m+1)>1$, we can expand the zeta function and reverse the order of summation and integration to obtain
\bee\label{e25}
I_L=(-)^m \pi^{2m-\fr} \sum_{n=1}^\infty n^{2m-1}\,K_n(a;m),
\ee
where
\[K_n(a;m):=\frac{1}{2\pi i}\int_L\,\frac{\g(s-m+\fs) \g(s-m+1)}{\g(s+1)}\,\left(\frac{a}{\pi^2n^2}\right)^s ds,\]
and we have employed the duplication formula for the gamma function
\[\g(2z)=2^{2z-1}\pi^{-\fr}\,\g(z) \g(z+\fs).\]

The integrals $K_n(a;m)$ have no poles in the half-plane $\Re (s)>m$, so that we can displace the path $L$ as far to the right as we please. On such a displaced path $|s|$ is everywhere large. The quotient of gamma functions may then be expanded by making use of the result given in \cite[p.~53]{PK}
\bee\label{e25a}
\frac{\g(s-m+\fs)\g(s-m+1)}{\g(s+1)}=\sum_{j=0}^{M-1} (-)^j c_j \g(s+\vartheta-j)+\rho_M(s) \g(s+\vartheta-M)
\ee
for positive integer $M$, where $\vartheta=\fs-2m$,
\[c_j=\frac{(m)_j (m+\fs)_j}{j!}=\frac{2^{-2j}(2m)_{2j} }{j!}\]
and $\rho_M(s)=O(1)$ as $|s|\ra\infty$ in $|\arg\,s|<\pi$. Substitution of this expansion into the integrals $K_n(a;m)$ then produces
\begin{eqnarray}
K_n(a;m)&=&\sum_{j=0}^{M-1} (-)^jc_j\,\frac{1}{2\pi i}\int_L \g(s+\vartheta-j)\,\left(\frac{a}{\pi^2n^2}\right)^{s}ds+{\cal R}_M\nonumber\\
&=&\sum_{j=0}^{M-1} (-)^j c_j \left(\frac{a}{\pi^2n^2}\right)^{2m+j-\fr} e^{-\pi^2n^2/a}+{\cal R}_M \label{e26}
\end{eqnarray}
by (\ref{e21}), where
\[{\cal R}_M=\frac{1}{2\pi i}\int_L\rho_M(s) \g(s+\vartheta-M) \left(\frac{a}{\pi^2n^2}\right)^sds.\]
Bounds for the remainder ${\cal R}_M$ have been considered in \cite[p.~71, Lemma 2.7]{PK}, where it is shown that
\bee\label{e27}
{\cal R}_M=O\left(\left(\frac{a}{\pi^2 n^2}\right)^{M-\vartheta} e^{-\pi^2 n^2/a} \right)
\ee
as $a\ra 0$ in the sector $|\arg\,a|<\fs\pi$.

Collecting together the results in (\ref{e24}), (\ref{e25}), (\ref{e26}) and (\ref{e27}), we obtain
\[I_L=(-)^m \left(\frac{a}{\pi}\right)^{2m-\fr} \sum_{n=1}^\infty \frac{e^{-\pi^2n^2/a}}{n^{2m}} \left\{\sum_{j=0}^{M-1}  c_j
\left(\!\frac{-a}{\ \pi^2n^2}\right)^j+O\left(\left(\frac{a}{\pi^2n^2}\right)^M\right)\right\}.\]
From (\ref{e23}) we now have the following theorem:
\newtheorem{theorem}{Theorem}
\begin{theorem}$\!\!\!.$
Let $m$ and $M$ be positive integers. Then, when $w=2m$, we have the expansion valid as $a\ra 0$ in $|\arg\,a|<\fs\pi$
\[S(a;2m)=\fs\g(\fs-m) a^{m-\fr}+\sum_{k=0}^m \frac{(-)^k}{k!}\,\zeta(2m-2k)\,a^k\hspace{4cm}\]
\bee\label{e28}
\hspace{6cm}+(-)^m\left(\frac{a}{\pi}\right)^{2m-\fr}\sum_{n=1}^\infty \frac{\Upsilon_n(a;m)}{n^{2m}}\,e^{-\pi^2n^2/a},
\ee
where $\Upsilon_n(a;m)$ has the asymptotic expansion
\[\Upsilon_n(a;m)=\sum_{j=0}^{M-1}\frac{(m)_j(m+\fs)_j}{j!} \,\left(\!\!\frac{-a}{\  \pi^2n^2}\right)^j+O\left(\left(\frac{a}{\pi^2n^2}\right)^M\right).
\]
\end{theorem}

This is the analogue of the Poisson-Jacobi transformation in (\ref{e11}). In the case $m=0$, the quotient of gamma functions in (\ref{e25a}) is replaced by the single gamma function $\g(s+\fs)$, with the result that $c_0=1$, $c_j=0$ ($j\geq 1$) and $\Upsilon_n(a;m)=1$ for all $n\geq 1$. Then (\ref{e28}) reduces to (\ref{e11}) and is valid for all values of the parameter $a$ (not just $a\ra 0$) satisfying $|\arg\,a|<\fs\pi$. 
\vspace{0.2cm}

\noindent{\bf Remark 1.}\ \ 
We note that the values of the zeta function appearing in (\ref{e28}) can be expressed alternatively in terms of Bernoulli numbers by the result \cite[p.~605]{DLMF}
\[\zeta(2n)=\frac{(2\pi)^{2n}}{2 (2n)!}\,|B_{2n}|.\]

\vspace{0.6cm}

\begin{center}
{\bf 4.\ Numerical results and concluding remarks}
\end{center}
\setcounter{section}{4}
\setcounter{equation}{0}
\renewcommand{\theequation}{\arabic{section}.\arabic{equation}}
From the well-known values \cite[p.~605]{DLMF}
\[\zeta(2)=\frac{\pi^2}{6},\qquad \zeta(4)=\frac{\pi^4}{90},\]
we obtain from Theorem 1 the expansions in the cases $m=1$ and $m=2$ given by
\bee\label{e31}
S(a;2)=\frac{\pi^2}{6}+\frac{a}{2}-(\pi a)^\fr-\left(\frac{a}{\pi}\right)^{\!\frac{3}{2}}\sum_{n=1}^\infty\frac{e^{-\pi^2n^2/a}}{n^2}
\bl\{\sum_{j=0}^{M-1}(\f{3}{2})_j \left(\!\frac{-a}{\ \pi^2n^2}\right)^j+O(a^M)\br\}
\ee
and
\[S(a;4)=\frac{\pi^4}{90}-\frac{\pi^2a}{6}-\frac{a^2}{4}+\frac{2}{3}\pi^\fr a^\frac{3}{2}\hspace{7cm}\]
\bee\label{e32}
\hspace{2cm}+\left(\frac{a}{\pi}\right)^{\!\frac{7}{2}}\sum_{n=1}^\infty\frac{e^{-\pi^2n^2/a}}{n^4}
\bl\{\sum_{j=0}^{M-1}\frac{(\f{5}{2})_j(2)_j}{j!} \left(\!\frac{-a}{\ \pi^2n^2}\right)^j+O(a^M)\br\}
\ee
valid as $a\ra 0$ in $|\arg\,a|<\fs\pi$.

\begin{table}[th]
\caption{\footnotesize{Values of the absolute error in the computation of $S(a;4)$ from (\ref{e32}).
The value of the index $j_0$ corresponds to optimal truncation of the expansion $\Upsilon_1(a;2)$.}}
\begin{center}
\begin{tabular}{|l|l|l|r|}
\hline
&&&\\[-0.25cm]
\mcol{1}{|c|}{$a$} & \mcol{1}{c|}{$S(a;4)$} & \mcol{1}{c|}{Error}  & \mcol{1}{c|}{$j_0$}\\
[.1cm]\hline
&&&\\[-0.3cm]
0.10 & 0.952696 & $9.662\times 10^{-86}$ &  96 \\
0.20 & 0.849025 & $9.768\times 10^{-43}$ &  46 \\
0.25 & 0.803169 & $4.045\times 10^{-34}$ &  36 \\
0.50 & 0.615128 & $7.769\times 10^{-17}$ &  17 \\
0.75 & 0.475493 & $4.656\times 10^{-11}$ &  10 \\
1.00 & 0.369026 & $3.642\times 10^{-8}$  &  6 \\
1.50 & 0.223285 & $2.856\times 10^{-5}$  &  3 \\
2.00 & 0.135356 & $7.500\times 10^{-4}$  &  1 \\
[.2cm]\hline
\end{tabular}
\end{center}
\end{table}
In Table 1 we show the results of numerical calculations for the case $m=2$. For different values of the parameter $a$ we present the value of the absolute error in the computation of $S(a;4)$ from (\ref{e32}). In the computations, we have used only the $n=1$ term (since the order of $\f{1}{4}\exp (-4\pi^2/a)$ was found to be less than the error), with the expansion for $\Upsilon_1(a;2)$ optimally truncated (corresponding to truncation at, or near, the least term in modulus) at index $j_0\simeq (\pi^2/a)-\f{5}{2}$.
It is seen that the error when $a=0.1$ is extremely small and that, only when $a\simeq 2$ does the relative error start to become significant.

To conclude, we mention that a similar treatment can be carried out for the sum
\[S_p(a;w)\equiv\sum_{n=1}^\infty \frac{e^{-an^p}}{n^{w}}\qquad (a\ra 0,\ \Re (a)>0)\]
for positive even integer $w$ and $p$. The case $w=0$ and $p>0$, corresponding to the Euler-Jacobi series, has been considered in \cite[\S 8.1]{PK}; see also \cite{K} for a hypergeometric approach when $p$ is a rational fraction. 
The details of the small-$a$ expansion of $S_p(a;w)$ will be presented elsewhere.

\vspace{0.6cm}

\noindent{\bf Acknowledgement.}\ \ The author wishes to acknowledge B. Guralnik for having brought this problem to his attention
\vspace{0.6cm}

\begin{center}
{\bf Appendix: \ A bound for the remainder $R_N$}
\end{center}
\setcounter{section}{1}
\setcounter{equation}{0}
\renewcommand{\theequation}{\Alph{section}.\arabic{equation}}
Let $\psi=\arg\,a$ and integer $N>\fs w+\fs$. Upon replacement of $s$ by $-s$ followed by use of (\ref{e23a}), the remainder $R_N$ in (\ref{e22a}) becomes
\[R_N=\frac{(2\pi)^w}{2\pi i}\int_{N-\fr-\infty i}^{N-\fr+\infty i} \zeta(1-w+2s)\,\frac{\g(1-w+2s)}{\g(1+s)}\,\frac{\sin \pi(s-\fs w)}{\sin \pi s}\left(\frac{a}{4\pi^2}\right)^{\!s}ds.\]
With $s=N-\fs+it$, $t\in(-\infty,\infty)$ we have
\[|R_N|\leq (2\pi)^{w-1} \left(\frac{a}{4\pi^2}\right)^{N-\fr} \zeta(2N-w) \int_{-\infty}^\infty
e^{-\psi t} \left|\frac{\g(2N-w+2it)}{\g(N+\fs+it)}\right|\,dt,\]
since $|\zeta(x+it)|\leq \zeta(x)$ ($x>1$) and
\[\left|\frac{\sin \pi(N-\fs-\fs w+it)}{\sin \pi(N+\fs+it)}\right|=\frac{|\cos \pi(\fs w-it)|}{\cosh \pi t}=
\frac{(\cos^2 \fs\pi w+\sinh^2 \pi t)^\fr}{\cosh \pi t}\leq 1.\]
It then follows that
\bee\label{a1}
|R_N|=O\left(\left(\frac{a}{\pi^2}\right)^{\!N-\fr} \int_{-\infty}^\infty e^{-\psi t}|\g(N-\fs w+it)|\,\left|\frac{\g(N-\fs w+\fs+it)}{\g(N+\fs+it)}\right|dt\right).
\ee

Using the argument presented in \cite[p.~126]{PK}, we set $N-\fs w-\fs=M+\delta$, with $-\fs<\delta\leq\fs$ so that $M\leq N-1$, to find
\[\left|\frac{\g(N-\fs w+\fs+it)}{\g(N+\fs+it)}\right|=P(t) g(t),\qquad g(t):=\left|\frac{\g(1+\delta+it)}{\g(\fs+it)}\right|,\]
where
\[P(t)=(\f{1}{4}+t^2)^{-\fr}\frac{\prod_{r=1}^M\{(r+\delta)^2+t^2\}^\fr}{\prod_{r=1}^{N-1}\{(r+\fs)^2+t^2\}^{\fr}}\leq
(\f{1}{4}+t^2)^{-\fr}\leq 2.\]
From the upper bound for the gamma function $\g(z)$ with $z=x+it$, $x>0$ \cite[p.~35]{PK}
\begin{eqnarray*}
|\g(z)|&<&\g(x) (1+t^2/x^2)^{\fr x-\frac{1}{4}}e^{-|t|\phi(t)}\, e^{1/(6|z|)},\qquad \phi(t)=\arctan (|t|/x)\\
\\
&<&\g(x) (1+\tau^2)^{\fr x-\frac{1}{4}} \exp\,[x\{\omega(\tau)-\fs\pi|\tau|\}\,e^{1/(6x)}\\
\\
&<&e^x \g(x) (1+\tau^2)^{\fr x-\frac{1}{4}} e^{-\fr\pi x|\tau|}\,e^{1/(6x)},
\end{eqnarray*}
where we have put $\tau=t/x$, defined $\omega(\tau)=|\tau| \arctan (1/|\tau|)$ and used the fact that
$0\leq\omega(\tau)<1$ for $\tau\in [0,\infty)$, with the limit 1 being approached as $\tau\ra\infty$.
Substituting the above bounds into (\ref{a1}), we see on setting $x=N-\fs w$ that
\[|R_N|=e^N \g(N\!-\!\fs w\!+\!1) \,O\left(\left(\frac{a}{\pi^2}\right)^{N-\fr}\int_0^\infty (1+\tau^2)^{N/2} g(\tau) \{e^{-\Delta_+\tau}+e^{-\Delta_-\tau}\}\,d\tau\right),\]
where $\Delta_\pm=(N-\fs w)(\fs\pi\pm\psi)$.
Since $g(\tau)=O(\tau^{\delta+\fr})$ as $\tau\ra\infty$, the integral is convergent provided $|\psi|<\fs\pi$ and is manifestly an increasing function of $N$.

Hence
\bee\label{a2}
R_N=O(a^{N-\fr})\qquad (a\ra 0,\ |\arg\,a|<\fs\pi),
\ee
with the constant implied in the $O$-symbol growing at least like $\g(N\!+\!1\!-\!\fs w)$ as $N$ increases.

\end{document}